\documentclass[10pt, twoside]{article}
\usepackage{amssymb,latexsym,amsmath}
\usepackage[hypertex]{hyperref}
\usepackage[english]{babel}
\date{}

\begin{document}
\makeatletter
\renewcommand{\@evenfoot}{ \thepage \hfil \footnotesize{\it 
ISSN 1025-6415 \ \ 
Dopovidi Natsionalno\" \i \ Akademi\" \i \  Nauk Ukra\"\i ni. 1997, no. 9} } 
\renewcommand{\@oddfoot}{\footnotesize{\it 
ISSN 1025-6415 \ \
Dopovidi Natsionalno\" \i \ Akademi\" \i \  Nauk Ukra\"\i ni. 1997, no. 9} 
\hfil \thepage } 
\noindent
\\
\\
\\
\\
\\
\\
\\
\\
\\
\\
\\
\\
\\
\par{
\leftskip=1.5cm  \rightskip=0cm  
\noindent 
UDC 512 \medskip \\
\copyright \ {\bf 1997}\medskip \\
{\bf T. R. Seifullin}\medskip \\ 
{\large \bf 
Homology of the Koszul complex of a system \smallskip \\
of polynomial equations \medskip \\ } 
{\it (Presented by Corresponding Member of the NAS of Ukraine A. A. Letichevsky)} \medskip \\
\small {\it 
Explicit complex morphism of a dual complex to the Koszul complex into the Koszul
complex is constructed. If the number of common roots of polynomials is finite 
in an algebraically closed field, then this mapping is a homotopic equivalence, 
thus explicit duality of the Koszul complex is obtained. }
\par \medskip } \noindent 
In the present paper it has generalized author's results of [1,2]  
on the whole  Koszul complex of polynomials 
of a system of polynomial equations,  
for the arbitrary number of polynomials, greater equal 
to  the number of variables.  In  result
it has constructed explicit complex morphism of a complex  dual to  
the Koszul complex into the Koszul complex. 
In the case of the finite number of common roots of polynomials in
algebraically  closed  field, this  mapping is a homotopic equivalence,  
thus it is obtained explicit duality for the Koszul complex.

{\bf  Koszul complexes.}
Let ${\bf R}$  be a commutative ring,  
   $f(x)=(f_1(x),\ldots  ,f_s(x))  \in {\bf R}[x]^s$   
be polynomials  in  
   $x=(x_1,\ldots  ,x_n)$  
with  coefficients  in a commutative ring 
   ${\bf R}$, $\forall i: f_i(x) \equiv   0$   
be a system of polynomial equations. 
Represent an element 
   $a=(a_1,\ldots  ,a_s)  \in {\bf R}^s$ 
in the form 
   $a\widehat{f}_*= \sum\limits_{i}^{}a_i\widehat{f}_*^i$, and $\forall  i$
denote by 
  $\widehat{f}_i$  
a linear  functional defined on 
  ${\bf R}^s$, 
assigning 
   $(a_1,\ldots  ,a_s)   \mapsto    a_i$,   
where   
   $\widehat{f}=(\widehat{f}_1,\ldots,\widehat{f}_s)$  
and  
   $\widehat{f}_*=(\widehat{f}_*^1,\ldots  ,\widehat{f}_*^s)$  be systems
of variables. Any functional 
is represented  in the form  
   $\widehat{f}b=  \sum\limits_{j}^{}\widehat{f}_jb^j$ with the action $a\widehat{f}_*.\widehat{f}b= $
   $ \sum\limits_{i}^{}a_i\widehat{f}_*^i.\sum\limits_{j}\widehat{f}_j b^j=  
   \sum\limits_{i}a_ib^i= ab$.
Then the element  
   $f(x)=(f_1(x),\ldots    ,f_s(x))$  
is represented in the form
   $f(x)\widehat{f}_*=  \sum\limits _{i}^{}f_i(x)\widehat{f}_*^i$.
Denote by 
   $(f(x))_x$ 
the ideal of the ring 
   ${\bf R}[x]$ 
generated by polynomials  
    $f(x)$.

For any graded module 
   ${\bf P}= \bigoplus\limits_r{\bf P}_r$,
denote  
   ${\bf P}_\times = \bigcup\limits_r{\bf P}_r$. 
If 
   $p\in {\bf P}_r$, 
then we will   write   
   $|p|=r$.   
Put   
   $\forall i:  |\widehat{f}_*^i|=-1$, $|\widehat{f}_i|=1$. 
Let  
   $\Lambda  (\widehat{f}_*)= \bigoplus\limits _r\Lambda _{-r}(\widehat{f}_*)$ 
be a    
Grassmann algebra, \hbox{i.e. an associative algebra}  with  $1$, 
free generated by elements\\ 
   $(\widehat{f}_*^1,\ldots ,\widehat{f}_*^s)$   
with relations  
   $\left\{ \left( \sum\limits _{i}^{}f_i\widehat{f}_*^i  \right) \left(  
   \sum\limits_{i}^{}f_i\widehat{f}_*^i \right) = 0\right\}$, 
where 
   $f=(f_1,\ldots ,f_s)$ 
are mutually commuting variables, 
   $\Lambda _{-r}(\widehat{f}_*)$ 
is the set of sums of products of $r$ elements in $(\widehat{f}_*)$.

Denote by 
   $\Lambda  _r(\widehat{f}_*)^*$  
the module dual to  
   $\Lambda _{-r}(\widehat{f}_*)$,  
i. e.  the set of  ${\bf R}$-linear  maps
   $\Lambda _{-r}(\widehat{f}_*) \rightarrow  {\bf R}$. 
Define a product in
   $\Lambda      (\widehat{f}_*)^*$ 
by \smallskip

\ \ \  $\forall c_1(\widehat{f}),c_2(\widehat{f})\in  
  \Lambda_{\times} (\widehat{f}_*)^*: 
  \forall a(\widehat{f}_*)\in
  \Lambda_{\times} (\widehat{f}_*): 
  a(\widehat{f}_*).c_1(\widehat{f})\cdot c_2(\widehat{f})=
  a(\widehat{f}'_*+\widehat{f}''_*).c_1(\widehat{f}')\cdot c_2(\widehat{f}''),$
\smallskip \\
where
  $\forall a_1(\widehat{f}_*),a_2(\widehat{f}_*)\in \Lambda  _{\times} (\widehat{f}_*): 
  a_2(\widehat{f}''_*)\cdot a_1(\widehat{f}'_*).c_1(\widehat{f}')\cdot
  c_2(\widehat{f}'')= \left( a_1(\widehat{f}_*).c_1(\widehat{f}) \right) \left(
  a_2(\widehat{f}_*).c_2(\widehat{f}) \right)$,
  \\
  $\widehat{f}' \simeq  \widehat{f}'' \simeq  \widehat{f}$ 
are equivalent systems of variables. 
Then $\Lambda _r(\widehat{f}_*)^* \simeq  \Lambda _r(\widehat{f})$.

Define
   ${\bf C}_{-r}(x,\widehat{f}_*)= {\bf R}[x]\otimes \Lambda _{-r}(\widehat{f}_*)$,\ \ 
   ${\bf C}(x,\widehat{f}_*)= \bigoplus\limits_r{\bf C}_{-r}(x,\widehat{f}_*)$,\ \ 
then
   $f(x)\widehat{f}_*=$ $\sum\limits _if_i(x)\widehat{f}_*^i $ $\in $
   ${\bf C}_{-1}(x,\widehat{f}_*)$.
Define
   ${\bf C}(x,\widehat{f})=    \bigoplus\limits _r{\bf C}_r(x,\widehat{f})$,
where     
   ${\bf C}_r(x,\widehat{f})$ 
is the set of ${\bf R}$-linear maps of  
   $\Lambda _{-r}(\widehat{f}_*)$  into  ${\bf R}[x]$,
which are equivalent to  
skew symmetric ${\bf R}$-polylinear forms  on 
   ${\bf R}^s \simeq  {\bf R}\otimes (\widehat{f}_*)$ into ${\bf R}[x]$.  
%There 
It holds  
   ${\bf C}_r(x,\widehat{f}) \simeq  {\bf R}[x]\otimes \Lambda _r(\widehat{f})$.

For $c(x,\widehat{f})\in   {\bf C}_r(x,\widehat{f})$,  if  $r\geq  1$  
denote by 
    $\partial \left[ c(x,\widehat{f})\right] $  
such  element
    $\in   {\bf C}_{r-1}(x,\widehat{f})$, 
that \ \ 
    $\forall a(x,\widehat{f}_*)\in {\bf C}_{-r+1}(x,\widehat{f}_*): \ \ 
    a(x,\widehat{f}_*).\partial \left[ c(x,\widehat{f})\right] =
    a(x,\widehat{f}_*)\cdot  f(x)\widehat{f}_*.c(x,\widehat{f})$,  
\ if \  $r=0$  \
denote 
    $\partial \left[ c(x,\widehat{f})\right]=0$.
Then 
    $\forall  c(x,\widehat{f})\in {\bf C}_r(x,\widehat{f}): 
    \partial \left[ \partial \left[  c(x,\widehat{f})\right] \right]     =0$.
Denote by \
    $({\bf C}(x,\widehat{f});\partial )= $
    $({\bf C}(x,\widehat{f});\widehat{f}\mapsto  f(x))$ 
a complex, \ which  is called the Koszul complex,  
    ${\bf Z}(x,\widehat{f})=\{ c\in  {\bf C}(x,\widehat{f})|$
    \hbox{$\partial
    \left[ c\right] =0\} $}, ${\bf B}(x,\widehat{f})=\{  \partial  \left[  c\right]
    |c\in  {\bf C}(x,\widehat{f})\}   $.   
From $\partial ^2=0$ it follows that
    ${\bf B}_r(x,\widehat{f})    \subseteq     {\bf Z}_r(x,\widehat{f})$.     
Denote by
    ${\bf H}_r(x,\widehat{f})= {\bf Z}_r(x,\widehat{f})/{\bf B}_r(x,\widehat{f})$.

%There 
It hold \smallskip

\ \ \ $ \forall  c_1,c_2\in {\bf C}_\times (x,\widehat{f}): 
    \partial \left[   c_1\cdot c_2\right] =
    \partial \left[ c_1\right] \cdot  c_2+  (-1)^{|c_1|}c_1\cdot
    \partial \left[ c_2\right],$

\ \ \ $ {\bf Z}(x,\widehat{f})\cdot  {\bf Z}(x,\widehat{f})\subseteq   {\bf Z}(x,\widehat{f}),
    \ {\bf B}(x,\widehat{f})\cdot  {\bf Z}(x,\widehat{f})\subseteq   {\bf B}(x,\widehat{f}),
    \ {\bf Z}(x,\widehat{f})\cdot {\bf B}(x,\widehat{f})\subseteq {\bf B}(x,\widehat{f}).
    $ \smallskip
 
Let 
    $v=(v_1,\ldots,v_m)$ and $\widehat{h}=(\widehat{h}_1,\ldots,\widehat{h}_t)$
be systems  of variables, 
    $h(v)=(h_1(v),\ldots   ,h_t(v))\in$ ${\bf R}[v]^t$, and
    $({\bf C}(v,\widehat{h});\ \partial )= ({\bf C}(v,\widehat{h});
     \widehat{h}\mapsto  h(v))$
be \ the \ Koszul \ complex \ 
of \ a \ system \ of \hbox{polynomials} $h(v)$. \   
\hbox{Denote} \  by \    
    $({\bf C}(v,\widehat{h},x_*,\widehat{f}_*); \partial ) \ =  
    ({\bf C}(v,\widehat{h},x_*,\widehat{f}_*); \ \widehat{f}\mapsto f(x), 
    \widehat{h}\mapsto h(v))$    
\ the \  complex \  of \   
    ${\bf R}$-linear \ maps \ ${\bf C}(x,\widehat{f})
    \rightarrow {\bf C}(v,\widehat{h})$,  
\ in \ which \ the \ boundary \ operator \ \hbox{defined} \  
as\ following:  
    \hbox{$\forall c\in {\bf C}_\times (x,\widehat{f}):  
    \forall a\in {\bf C}_\times (v,\widehat{h},x_*,\widehat{f}_*):$}     
    $\partial \left[ a\right] .c = \partial
    \left[ a.c\right] - (-1)^{|a|}a.\partial \left[ c\right]$,
where 
    ${\bf C}_r(v,\widehat{h},x_*,\widehat{f}_*)= 
    \{ a\in {\bf C}(v,\widehat{h},x_*,\widehat{f}_*)|  
    \forall r': a.{\bf C}_{r'}(x,\widehat{f}) 
    \subseteq {\bf C}_{r'+r}(v,\widehat{h})\}$.  
Denote by 
    ${\bf Z}(v,\widehat{h},x_*,\widehat{f}_*)=$
    \hbox{$\{ a\in {\bf C}(v,\widehat{h},x_*,\widehat{f}_*)|$}
    $\partial \left[ a\right] =0\},
    {\bf B}(v,\widehat{h},x_*,\widehat{f}_*)= 
    \{ \partial \left[ a\right] |a\in $
    ${\bf C}(v,\widehat{h},x_*,\widehat{f}_*)\}$,
    ${\bf H}(v,\widehat{h},x_*,\widehat{f}_*)=$
    ${\bf Z}(v,\widehat{h},x_*,\widehat{f}_*)/ 
    {\bf B}(v,\widehat{h},x_*,\widehat{f}_*)$.

Elements of 
    ${\bf Z}(v,\widehat{h},x_*,\widehat{f}_*)$ 
are called  complex morphisms,
elements of 
    ${\bf B}(v,\widehat{h},x_*,\widehat{f}_*)$ 
are maps homotopic to zero.

If 
    $v=()$, $h(v)=()$ 
and 
    $\widehat{h}=()$, 
then
    $({\bf C}(x_*,\widehat{f}_*);\partial ) =$
    $({\bf C}(v,\widehat{h},x_*,\widehat{f}_*);\partial )$ 
is the complex dual to the Koszul complex 
    $({\bf C}(x,\widehat{f});\partial )$.
If 
    $x=()$, $f(x)= ()$ 
and 
    $\widehat{f}= ()$, 
then
    $({\bf C}(v,\widehat{h},x_*,\widehat{f}_*);\partial)= $
$ ({\bf C}(v,\widehat{h});\partial )$.

Let 
    $({\bf C}(w,\widehat{g},v_*,\widehat{h}_*);\partial )=$
    $({\bf C}(w,\widehat{g},v_*,\widehat{h}_*)$; $\widehat{g}\mapsto  g(w),  
    \widehat{h}\mapsto h(v))$
be a complex of  ${\bf R}$-linear  maps 
    ${\bf C}(v,\widehat{h})\rightarrow  {\bf C}(w,\widehat{g})$, 
then \smallskip

\ \ \ $ \forall a\in {\bf C}_\times (w,\widehat{g},v_*,\widehat{h}_*): 
    \forall b\in {\bf C}_\times (v,\widehat{h},x_*,\widehat{f}_*): 
    \partial \left[ a.b \right] =
    \partial \left[ a \right] .b+ (-1)^{|a|}a.\partial \left[ b\right],$
  
\ \ \ ${\bf Z}(w,\widehat{g},v_*,\widehat{h}_*).{\bf Z}(v,\widehat{h},x_*,\widehat{f}_*)
    \subseteq  {\bf Z}(w,\widehat{g},x_*,\widehat{f}_*),$
\eject

\ \ \ ${\bf B}(w,\widehat{g},v_*,\widehat{h}_*).{\bf Z}(v,\widehat{h},x_*,\widehat{f}_*)
    \subseteq  {\bf B}(w,\widehat{g},x_*,\widehat{f}_*),$
    
\ \ \ ${\bf Z}(w,\widehat{g},v_*,\widehat{h}_*).{\bf B}(v,\widehat{h},x_*,\widehat{f}_*)
    \subseteq  {\bf B}(w,\widehat{g},x_*,\widehat{f}_*).$ \smallskip

Let 
    $a(v,\widehat{h}),b(v,\widehat{h})\in{\bf C}_\times(v,\widehat{h})$ 
and
    $(-1)^{|a|}=1$,$(-1)^{|b|} = -1$,
denote by 
    ${\bf 1}_{(x,\widehat{f})}(a(v,\widehat{h}),b(v,\widehat{h}))$
the map of the form \smallskip

\ \ \ ${\bf C}(x,\widehat{f}) \ni  c(x,\widehat{f}) \mapsto
    {\bf 1}_{(x,\widehat{f})}(a(v,\widehat{h}),b(v,\widehat{h})).c(x,\widehat{f})=
    c(a(v,\widehat{h}),b(v,\widehat{h})) \in  {\bf C}(v,\widehat{h}).$\smallskip

Define the product \smallskip

\ \ \ ${\bf C}_\times (x,\widehat{f},v_*,\widehat{h}_*)\times {\bf C}_\times (x,\widehat{f},w_*,\widehat{g}_*)
   \ni (c,c') \mapsto c\cdot c' \in
   {\bf C}_\times (x,\widehat{f},v_*,\widehat{h}_*,w_*,\widehat{g}_*): $

\ \ \ $\forall a\in {\bf C}_\times (v,\widehat{h}): \forall a'\in {\bf C}_\times (w,\widehat{g}): 
   c\cdot c'.a\cdot a' = (-1)^{|a||c'|}(c.a)\cdot (c'.a').$ 
\smallskip \\
Then \smallskip

\ \ \ $\partial \left[ c\cdot c'\right] = \partial \left[  c\right]
   \cdot c'+ (-1)^{|c|}c\cdot \partial \left[ c'\right].$ \smallskip

Let  
   ${\bf C}_\times (x,\widehat{f},x_*,\widehat{f}_*,w,\widehat{g},v_*,\widehat{h}_*)  \ni
   c(x,\widehat{f},x_*,\widehat{f}_*,w,\widehat{g},v_*,\widehat{h}_*)=$\\
\indent \hphantom{Let} 
   $p(w,\widehat{g},v_*,\widehat{h}_*)\cdot a(x_*,\widehat{f}_*)\cdot a'(x,\widehat{f}) \in
   {\bf C}_\times (w,\widehat{g},v_*,\widehat{h}_*)\cdot {\bf C}_\times (x_*,\widehat{f}_*)\cdot
   {\bf C}_\times (x,\widehat{f}),$
   \\
denote by
   $\mathop{\top}\limits_{(x,\widehat{f})}  c(x,\widehat{f},x_*,\widehat{f}_*,w,\widehat{g},v_*,\widehat{h}_*)=
   p(w,\widehat{g},v_*,\widehat{h}_*)\cdot   \left(   a(x_*,\widehat{f}_*).a'(x,\widehat{f})
   \right),$
   \\
denote by
   $ \mathop{\bot}\limits_{(x,\widehat{f})} c(x,\widehat{f},x_*,\widehat{f}_*,w,\widehat{g},v_*,\widehat{h}_*)$
element of
   ${\bf C}_\times (x_*,\widehat{f}_*,w,\widehat{g},v_*,\widehat{h}_*)$
such  that
   $\forall b(x,\widehat{f})\in {\bf C}_\times (x,\widehat{f}): $
   $\mathop{\bot}\limits_{(x,\widehat{f})}
   c(x,\widehat{f},x_*,\widehat{f}_*,w,\widehat{g},v_*,\widehat{h}_*).b(x,\widehat{f})=
   p(w,\widehat{g},v_*,\widehat{h}_*)\cdot
   \left( a(x_*,\widehat{f}_*).a'(x,\widehat{f})b(x,\widehat{f}) \right) $

An exponential determinant we call \smallskip

\ \ \ $\det   \left\| \begin{matrix} \ bc& B\widehat{f}_*\ \\ \ \widehat{h} C&  0 \  \end{matrix}
   \right\| =
   \mathop{\top}\limits_{\widehat{p}} (B^l\widehat{f}_*+b^l\widehat{p}_*)\cdot \ldots   \cdot
   (B^1\widehat{f}_*+b^1\widehat{p}_*)\cdot   (\widehat{h}C_1+\widehat{p}c_1)\cdot    \ldots
   \cdot (\widehat{h}C_m+\widehat{p}c_m),$ \smallskip \\
where    
   $\forall    k:   B^k\widehat{f}_*+b^k\widehat{p}_*\in     \Lambda     _{-
   1}(\widehat{f}_*,\widehat{p}_*)$, $\forall k: \widehat{h}C_k+\widehat{p}c_k\in   \Lambda
   _1(\widehat{h},\widehat{p})$, $\widehat{f},  \widehat{h},  \widehat{p}$ 
are collections of anticommuting variables.

{\bf Difference Jacobian.}
Let 
   $f(x)=(f_1(x),\ldots ,f_s(x)) \in {\bf R}\left[ x\right] ^s$
be polynomials in  
   $x=(x_1,\ldots ,x_n)$ with coefficients in  
a commutative ring  
   ${\bf R}$.  
Let $y \simeq x$, 
   \hbox{ $x-y= (x_1-y_1,\ldots ,x_n-y_n)$}
   $\in {\bf R}\left[ x,y\right] ^n$. 
Consider the complex  
   \hbox{$({\bf C}(x,y,\widehat{f}_x,\widehat{f}_y,\widehat{u});\partial)= 
   ({\bf C}(x,y,\widehat{f}_x,\widehat{f}_y,\widehat{u}); 
   \widehat{f}_x{\mapsto} f(x), \widehat{f}_y{\mapsto} f(y),$}
    $\widehat{u}{\mapsto} (x-y))$. 
%To abridge notations we will write 
To shorten notations we will write
   $p_x= (x,\widehat{f}_x)$ and $p_*^x= (x_*,\widehat{f}_*^x)$.

Define    
   $\triangle _{(x,\widehat{u})}(x,y,\widehat{u}).c(x,\widehat{u}) =
   c(x,\widehat{u}) - c(y,\widehat{0})$.

{\bf Lemma.} {\it There exists 
   $\nabla  _{(x',\widehat{u}')}(x,y,\widehat{u}) \in
   {\bf C}(x'_*,\widehat{u}'_*,x,y,\widehat{u})$ 
such that \smallskip

\ \ \ $\triangle_{(x',\widehat{u}')}(x,y,\widehat{u})  =   \partial   \left[   \nabla
   _{(x',\widehat{u}')}(x,y,\widehat{u})\right]               \in
   {\bf B}(x'_*,\widehat{u}'_*,x,y,\widehat{u}),$ 
\smallskip \\
the operator  
   $\nabla  _{(x',\widehat{u}')}(x,y,\widehat{u})$  
is called a difference homotopy operator. }

{\bf Proof.}  Set
   \begin{align*}
   \triangle ^k_{(x,\widehat{u})}(x,y,\widehat{u}).c(x,\widehat{u})=  
   &\hphantom{-}\ \,c(y_1,\widehat{0},\ldots ,y_{k-1},\widehat{0},x_k,\widehat{u}_k, x_{k+1},
   \widehat{u}_{k+1},\ldots ,x_n,\widehat{u}_n)-\\  
   &-c(y_1,\widehat{0},\ldots ,y_{k-1},\widehat{0},y_k,{\ }
   \widehat{0}_{\ },x_{k+1},\widehat{u}_{k+1},
   \ldots ,x_n,\widehat{u}_n).\hphantom{cccccccccc}
   \end{align*}
Then  
   $\triangle  =\sum\limits_{k}^{}\triangle  ^k=$
   $ \sum\limits_{k}\partial   \left[\frac{\widehat{u}_k}{x_k-y_k}\right] 
   \cdot \triangle  ^k=
   \sum\limits_{k}\partial \left[  \frac{\widehat{u}_k}{x_k-y_k}\cdot
   \triangle  ^k\right] 
   =   \partial
   \left[ \sum\limits_{k}\frac{\widehat{u}_k}{x_k-y_k}\cdot \triangle ^k\right]
   $, 
since
   \eject \noindent
   $\partial \left[ \triangle ^k\right] =0$.

%There 
It hold   
   $\partial \left[ \triangle \right]  = 0$, \ 
and  
   $\triangle
   .\triangle =\triangle  $.  
Then  
   $\partial  \left[  \triangle  .\nabla
   \right]  =  \triangle  .\partial  \left[  \nabla  \right]   =   \triangle
   .\triangle  =\triangle $, hence, if $\nabla $ 
is a difference homotopy, then 
   $\triangle .\nabla $ 
is also a difference homotopy,  moreover,
   $\triangle .(\triangle .\nabla )= \triangle .\nabla $. 
A difference homotopy 
   $\nabla$, for which   
%we have 
   $\triangle .\nabla = \nabla $, 
we call reduced. 
Let 
   $\nabla '$ and  $\nabla  ''$  
be two reduced difference homotopy, then 
   $\partial \left[ \nabla '-\nabla ''\right] =
   \triangle  -\triangle  =  0$,  hence,  $\nabla   '-\nabla   ''\in
   {\bf Z}$. Further $\nabla '-\nabla  ''=  \triangle  .(\nabla  '-\nabla  '')=
   \partial \left[ \nabla \right] .(\nabla  '-\nabla  '')=  \partial  \left[
   \nabla .(\nabla '-\nabla '')\right]  \in {\bf B}$,  
hence, a reduced difference homotopy is uniquely determined, up to homotopy. If
   $c(x,\widehat{u}) \in  {\bf Z}(x,\widehat{u})$, 
and 
   $\nabla '$, $\nabla ''$
are two reduced difference homotopy, then 
   $\nabla '.c- \nabla ''.c= \partial
   \left[ \nabla .(\nabla  '-\nabla  '')\right]  .c=\partial  \left[  \nabla
   .(\nabla '-\nabla '').c\right]  \in {\bf B}$,  
i. e.  a reduced difference homotopy of an element 
   $\in {\bf Z}(x,\widehat{u})$  
is uniquely determined, up to addend in 
   ${\bf B}(x,y,\widehat{u})$.

For a polynomial 
   $F(x) \in  {\bf R}[x]= $
   ${\bf C}_0(x,\widehat{u})$= ${\bf Z}_0(x,\widehat{u})$
denote by
   $\widehat{u}\nabla  F(x,y) =$ $\sum\limits_{k}^{} \widehat{u}_k\nabla ^kF(x,y)=$
   $\nabla_{(x,\widehat{u})}(x,y,\widehat{u}).F(x)$, \ 
then  the element  \ 
   $ \widehat{u}\nabla F(x,y)$ \ 
is a reduced \ difference \ homotopy of $F(x)$, \ since  
   $\triangle_{(x,\widehat{u})}(x,y,\widehat{u}).\widehat{u}\nabla F(x,y)= 
   \widehat{u}\nabla  F(x,y)-\widehat{0}\nabla F(y,y)=$ $\widehat{u}\nabla 
   F(x,y)$,
hence,
   $\widehat{u}\nabla F(x,y)$ 
is uniquely determined in ${\bf H}(x,y,\widehat{u})$.

A difference Jacobian is called \smallskip

\ \ \ $J(p_x,p_y)=
   \det \left\| \begin{matrix}\ \nabla f(x,y) \ \\ 
   \ \widehat{f}_x-\widehat{f}_y \ \end{matrix} \right\| =
   \mathop{\top}\limits_{\widehat{u}}  
   \det \| -\widehat{u}_*\|   \det   \|   \widehat{f}_x-\widehat{f}_y-
   \widehat{u}\nabla f(x,y)\|. $ \smallskip

{\bf Lemma.} {\it \smallskip

1. 
   $\partial \left[
   \det \left\| \begin{matrix}\ \nabla f(x,y)\ \\ 
   \ \widehat{f}_x-\widehat{f}_y \ \end{matrix} \right\| \right] =0$, i. e.
   $\det \left\| \begin{matrix}\ \nabla f(x,y)\ \\ 
   \ \widehat{f}_x-\widehat{f}_y \ \end{matrix} \right\|
   \in  {\bf Z}(x,y,\widehat{f}_x,\widehat{f}_y).$

2.     
A difference Jacobian is uniquely determined in
   ${\bf H}(x,y,\widehat{f}_x,\widehat{f}_y)$, 
independently of the choice of 
   $\nabla f(x,y)$.}

{\bf  Proof 1.} \smallskip

\ \ \ $\partial  \left[  \det  \|  -\widehat{u}_*\| \right] =
   -\sum\limits_{k}^{}(x_k-y_k)\widehat{u}_*^k\cdot \det \| -\widehat{u}_*\|  =0,$

\ \ \ $\partial  \left[  \det  \|  \widehat{f}_x-\widehat{f}_y-\widehat{u}\nabla   f(x,y)\|
   \right]   =    \partial    \left[    \prod_i(\widehat{f}_{i,x}-\widehat{f}_{i,y}-
   \widehat{u}\nabla f_i(x,y))\right] = 0,$ since

\ \ \ $\forall  i: \partial  \left[  \widehat{f}_{i,x}-\widehat{f}_{i,y}-\widehat{u}\nabla
   f_i(x,y)\right] = f_i(x)- f_i(y)-(x-y)\nabla f_i(x,y)= 0.$ 
\smallskip \\
Hence,
   $\partial \left[ \det \left\| \begin{matrix}\ \nabla f(x,y) \ \\ 
   \ \widehat{f}_x-\widehat{f}_y \ \end{matrix} \right\| \right] = 
   \partial \left[ \mathop{\top}\limits_{\widehat{u}} 
   \det \|  -\widehat{u}_*\|  \det  \|  \widehat{f}_x-\widehat{f}_y-\widehat{u}\nabla f(x,y)\| 
   \right] = 0.$

{\bf Proof 2.} Since 
   \ $\forall i: \ \widehat{u}\nabla f_i(x,y)$  
   \ \ is uniquely determined \ in \ ${\bf H}$, \ then  $\widehat{f}_{i,x}-\widehat{f}_{i,y}-
   \widehat{u}\nabla f_i(x,y)$ 
is uniquely determined in ${\bf H}$. 
Hence, \\
   $\det \left\| \begin{matrix} \ \nabla f(x,y) \ 
   \\ \ \widehat{f}_x-\widehat{f}_y \ \end{matrix} \right\| 
   = \mathop{\top}\limits_{\widehat{u}} \det \| -\widehat{u}_*\| \prod_i(\widehat{f}_{i,x}-\widehat{f}_{i,y}-
   \widehat{u}\nabla f_i(x,y))$ is uniquely determined  in  ${\bf H}$,  since
   $\forall  i:  \partial  \left[  \widehat{f}_{i,x}-\widehat{f}_{i,y}-\widehat{u}\nabla
   f_i(x,y)\right] =0$ and $\partial \left[ \det \| -\widehat{u}_*\| \right] = 0$.
   
Define
   $\det \left\| \begin{matrix} \ \nabla f(x,y) & -\widehat{u}_* \ \\ 
   \ \widehat{f}_x-\widehat{f}_y & 0 \ \end{matrix} \right\| = $
   $\mathop{\top}\limits_{\widehat{u}'}  \det  \|  -\widehat{u}_*-\widehat{u}'_*\|  \det   \|
   \widehat{f}_x-
   \widehat{f}_y-\widehat{u}'\nabla f(x,y)\| $. Then \smallskip

\ \ \ $\partial \left[
   \det \left\| \begin{matrix} \ \nabla f(x,y) & -\widehat{u}_*\ \\ 
   \ \widehat{f}_x-\widehat{f}_y & 0 \ \end{matrix} \right\|
   \right] =
   \partial \left[ \mathop{\top}\limits_{\widehat{u}'} \det \|  -\widehat{u}_*-\widehat{u}_*'\|  \det  \|
   \widehat{f}_x-\widehat{f}_y-\widehat{u}'\bigtriangledown f(x,y)\| \right] = 0.$
\eject

{\bf Theorem.}\smallskip 

\ \ \ $\det \left\| \begin{matrix} \ \nabla f(x,y) \ \\
   \widehat{f}_x-\widehat{f}_y \ \end{matrix} \right\| \left( {\bf 1}_{(z,\widehat{f}_z)}(x,\widehat{f}_x)-
   {\bf 1}_{(z,\widehat{f}_z)}(y,\widehat{f}_y) \right)=$

\ \ \ $=\ \ \,\partial \left[  \mathop{\top}\limits_{\widehat{x}'} \mathop{\top}\limits_{\widehat{u}'} \mathop{\top}\limits_{\widehat{u}}
   \det \left\| \begin{matrix}\ \nabla f(x,y) & \widehat{u}_*\ \\ 
   \ -\widehat{f}_x+\widehat{f}_y & 0 \ \end{matrix} \right\|
   \nabla _{(x',\widehat{u}')}(x,y,\widehat{u})\cdot
   {\bf 1}_{(z,\widehat{f}_z)}(x',\widehat{f}_y+\widehat{u}'\nabla f(x',y)) \right].$\smallskip 

{\bf Proof.} \smallskip

\ \ \ $\det \left\| \begin{matrix} \ \nabla f(x,y) \ \\ 
   \ \widehat{f}_x-\widehat{f}_y \ \end{matrix} \right\|
   \left({\bf 1}_{(z,\widehat{f}_z)}(x,\widehat{f}_x )-
   {\bf 1}_{(z,\widehat{f}_z)}(y,\widehat{f}_y) \right) =$

\ \ \ $\qquad = \mathop{\top}\limits_{\widehat{u}}  
   \det \| -\widehat{u}_*\| \det \|
   \widehat{f}_x-\widehat{f}_y-\widehat{u}\nabla   f(x,y)\|   
   \left( {\bf 1}_{(z,\widehat{f}_z)}(x,\widehat{f}_x)-
   {\bf 1}_{(z,\widehat{f}_z)}(y,\widehat{f}_y) \right) =$
 
\ \ \ $\qquad = \mathop{\top}\limits_{\widehat{u}}  
   \det  \|  -\widehat{u}_*\|   
   \det \| \widehat{f}_x-\widehat{f}_y-\widehat{u}\nabla f(x,y)\|
   \left( {\bf 1}_{(z,\widehat{f}_z)}(x,\widehat{f}_y+\widehat{u}\nabla 
   f(x,y))
   -{\bf 1}_{(z,\widehat{f}_z)}(y,\widehat{f}_y) \right) =$

\ \ \ $\quad \qquad = \mathop{\top}\limits_{\widehat{u}'} 
   \mathop{\top}\limits_{\widehat{u}} 
   \det \|  -\widehat{u}_*-\widehat{u}'_*\|  
   \det \| \widehat{f}_x-\widehat{f}_y-\widehat{u}'\nabla f(x,y)\| \cdot $

\ \ \ $\qquad \qquad \qquad \qquad \qquad \qquad \qquad  \qquad \qquad 
   \cdot \left( {\bf 1}_{(z,\widehat{f}_z)}(x,\widehat{f}_y+\widehat{u}\nabla f(x,y))-
   {\bf 1}_{(z,\widehat{f}_z)}(y,\widehat{f}_y) \right) =$

\ \ \ $\qquad = \mathop{\top}\limits_{\widehat{u}}
   \det \left\| \begin{matrix} \ \nabla f(x,y) & -\widehat{u}_* \ \\ 
   \ \widehat{f}_x-\widehat{f}_y & 0 \ \end{matrix} 
   \right\| \left( {\bf 1}_{(z,\widehat{f}_z)}(x,\widehat{f}_y+\widehat{u}\nabla f(x,y))-
   {\bf 1}_{(z,\widehat{f}_z)}(y,\widehat{f}_y) \right) =$

\ \ \ $\qquad = \mathop{\top}\limits_{\widehat{x}'} 
   \mathop{\top}\limits_{\widehat{u}'} \mathop{\top}\limits_{\widehat{u}} 
   \det \left\| \begin{matrix} \ \nabla f(x,y) & -\widehat{u}_*\ \\ 
   \ \widehat{f}_x-\widehat{f}_y & 0 \ \end{matrix} \right\| 
   \partial  \left[  \nabla  _{(x',\widehat{u}')}(x,y,\widehat{u})  \right]    \cdot
   {\bf 1}_{(z,\widehat{f}_z)}(x',\widehat{f}_y+\widehat{u}'\nabla f(x',y))=$

\ \ \ $\qquad = \partial \left[ \mathop{\top}\limits_{\widehat{x}'} 
   \mathop{\top}\limits_{\widehat{u}'} \mathop{\top}\limits_{\widehat{u}}
   \det \left\| \begin{matrix}\ \nabla f(x,y) & \widehat{u}_* \ \\ 
   \ -\widehat{f}_x+\widehat{f}_y & 0 \ \end{matrix} \right\|
   \nabla _{(x',\widehat{u}')}(x,y,\widehat{u})\cdot
   {\bf 1}_{(z,\widehat{f}_z)}(x',\widehat{f}_y+\widehat{u}'\nabla f(x',y)) \right].$
\smallskip 

Denote \smallskip

\ \ \ $T(p_*^z,p_x,p_y)=
   \mathop{\top}\limits_{\widehat{x}'} \mathop{\top}\limits_{\widehat{u}'} \mathop{\top}\limits_{\widehat{u}}
   \det \left\| \begin{matrix} \ \nabla f(x,y) & \widehat{u}_* \ \\ 
   \ -\widehat{f}_x+\widehat{f}_y & 0 \  \end{matrix} \right\|
   \nabla              _{(x',\widehat{u}')}(x,y,\widehat{u})\cdot
   {\bf 1}_{(z,\widehat{f}_z)}(x',\widehat{f}_y+\widehat{u}'\nabla f(x',y)). $
\smallskip

{\bf  Duality.}

{\bf Theorem-definition 1.} 
{\it $J$-map we call a map \smallskip

\ \ \ ${\bf C}(x_*,\widehat{f}_*^x) \ni c(x_* ,\widehat{f}_*^x) \mapsto
   \mathop{\top}\limits_{(y,\widehat{f}_y)}
   \det \left\| \begin{matrix}\ \nabla f(x,y)\ \\ 
   \ \widehat{f}_x-\widehat{f}_y \ \end{matrix} \right\|
   c(y_*,\widehat{f}_*^y) \in  {\bf C}(x,\widehat{f}_x)$ \smallskip

1. 
$J$-map is a complex morphism: \smallskip

\ \ \ $ \mathop{\top}\limits_{p_y}  \partial 
   \left[ J(p_x,p_y) \right] c(p_*^y)=0.$ \smallskip

2. 
$J$-map is uniquely determined,  up to  homotopy,
independently of the choice of $\nabla f(x,y)$: \smallskip

\ \ \ $\mathop{\top}\limits_{p_y} J(p_x,p_y)c(p_*^y)-\mathop{\top}\limits_{p_y} 
   J'(p_x,p_y)c(p_*^y) = \mathop{\top}\limits_{p_y} 
   \partial \left[ S(p_x,p_y) \right] c(p_*^y).$ \smallskip

3. 
$J$-map is a homotopy ${\bf C}(x,\widehat{f}_x)$-skewlinear: \smallskip 

\ \ \ $\left(   \mathop{\top}\limits_{p_y} J(p_x,p_y)c(p_*^y) \right) a(p_x)-
   \mathop{\top}\limits_{p_y} J(p_x,p_y) \left( \mathop{\bot}\limits_{p_y} c(p_*^y)a(p_y) \right) =
   \mathop{\top}\limits_{p_y}   \mathop{\top}\limits_{p_z}   \partial   \left[   T(p_*^z,p_x,p_y) \right]
   c(p_*^y)a(p_z).$ \smallskip }\eject

{\bf Theorem-definition 2.} {\it $J$-product is called a map \smallskip
 
\ \ \ ${\bf C}(x_*,\widehat{f}_*^x)\times {\bf C}(y_*,\widehat{f}_*^y) \ni \left(
   c_1(x_*,\widehat{f}_*^x),c_2(y_* ,\widehat{f}_*^y) \right) \mapsto$

\ \ \ $\mathop{\bot}\limits_{(x,\widehat{f}_x)} c_1(x_*,\widehat{f}_*^x) \mathop{\top}\limits_{(y,\widehat{f}_y)}
   \det \left\| \begin{matrix} \ \nabla f(x,y) \ \\ 
   \ \widehat{f}_x-\widehat{f}_y \ \end{matrix} \right\|
   c_2(y_*,\widehat{f}_*^y) \in  {\bf C}(x_*,\widehat{f}_*^x).$\smallskip

1. 
$J$-product is a complex bimorphism: \smallskip 

\ \ \ $\mathop{\bot}\limits_{p_x} c_1(p_*^x)\mathop{\top}\limits_{p_y}   \partial  \left[  J(p_x,p_y)  \right]
   c_2(p_*^y)= 0.$\smallskip 

2. 
$J$-product is uniquely determined  up to  homotopy,  
independently of the choice of 
\smallskip   $\nabla f(x,y)$: 

\ \ \ $\mathop{\bot}\limits_{p_x}   c_1(p_*^x)\mathop{\top}\limits_{p_y}   J(p_x,p_y)c_2(p_*^y)-   \mathop{\bot}\limits_{p_x}
   c_1(p_*^x)\mathop{\top}\limits_{p_y} J'(p_x,p_y)c_2(p_*^y)=
   \mathop{\bot}\limits_{p_x} c_1(p_*^x)  \mathop{\top}\limits_{p_y}  \partial  \left[  S(p_x,p_y)  \right]
   c_2(p_*^y).$ \smallskip

3. $J$-product is homotopy skewcommutative:\smallskip

\ \ \ $\mathop{\bot}\limits_{p_x} c_1(p_*^x)
   \mathop{\top}\limits_{p_y} J(p_x,p_y)c_2(p_*^y)-(-1)^{|c_1|'|c_2|'}
   \mathop{\bot}\limits_{p_x} c_2(p_*^x)
   \mathop{\top}\limits_{p_y} J(p_x,p_y)c_1(p_*^y)=$

\ \ \ \ \ \ $\quad =\mathop{\top}\limits_{p_y} 
   \mathop{\top}\limits_{p_z} c_1(p_*^y)\partial \left[ C(p_*^x,p_y,p_z) 
   \right]c_2(p_*^z),$\smallskip \\
where 
   $|c_1|'=|c_1|+|J|$ and $|c_2|'=|c_2|+|J|$.}

{\bf Theorem 3.} 
{\it If the image of ${\bf H}(x_*,\widehat{f}_*^x)$ 
under the  map $J$ include $1\in  
    {\bf H}(x,\widehat{f}_x)$, 
i. e. if \smallskip

\ \ \ $\exists e(x_*,\widehat{f}_*^x)\in {\bf Z}(x_*,\widehat{f}_*^x): 
    \mathop{\top}\limits_{(y,\widehat{f}_y)} 
    \det \left\| \begin{matrix}\ \nabla f(x,y) \ \\ 
    \ \widehat{f}_x-\widehat{f}_y \ \end{matrix} \right\|
    e(y_*,\widehat{f}_*^y)= 1+\partial \left[ t(x,\widehat{f}_x) \right],$
\smallskip \\
then the  map \smallskip

\ \ \ ${\bf C}(x,\widehat{f}_x)\ni    c(x,\widehat{f}_x)   \mapsto
    \mathop{\bot}\limits_{(x,\widehat{f}_x)} e(x_*,\widehat{f}_*^x)\,c(x,\widehat{f}_x)\in
    {\bf C}(x_* ,\widehat{f}_*^x)$
\smallskip \\
is homotopy inverse to the map $J$, i. e.

1)
   $\forall c(p_x)\in {\bf C}(p_x): c(p_x)-  \mathop{\top}\limits_{p_y}  J(p_x,p_y)  \left(
   \mathop{\bot}\limits_{p_y}  e(p_*^y)c(p_y)  \right)   =
   \mathop{\top}\limits_{p_y}   \partial   \left[ R(p_*^y,p_x) \right] c(p_y),$

2)
   $\forall c(p_*^x)\in {\bf C}(p_*^x): c(p_*^x)-\mathop{\bot}\limits_{p_x} e(p_*^x) \left(
   \mathop{\top}\limits_{p_y}  J(p_x,p_y)c(p_*^y)  \right)  =  \mathop{\top}\limits_{p_y}  \partial   \left[
   L(p_*^x,p_y) \right] c(p_*^y),$
   \\
and inverse to the map  $J$ in 
   ${\bf H}(x,\widehat{f}_x ) \rightarrow {\bf H}(x_*,\widehat{f}_*^x)$, i. e.

3)
   $\forall  c(p_x)\in  {\bf Z}(p_x): c(p_x)-  \mathop{\top}\limits_{p_y}   J(p_x,p_y)\left(
   \mathop{\bot}\limits_{p_y}  e(p_*^y)c(p_y)  \right)   =   \partial   \left[   \mathop{\top}\limits_{p_y}
   R(p_*^y,p_x)c(p_y) \right]$,

4)
   $\forall c(p_*^x)\in {\bf Z}(p_*^x): c(p_*^x)-\mathop{\bot}\limits_{p_x} 
   e(p_*^x)\left(
   \mathop{\top}\limits_{p_y}  J(p_x,p_y)c(p_*^y)  \right)   =\partial   \left[   
   \mathop{\top}\limits_{p_y}
   L(p_*^x,p_y )c(p_*^y) \right]$;
\\

5) 
   ${\bf H}_0(x,\widehat{f}_x)$   
and 
   ${\bf H}_{-s+n}(x_*,\widehat{f}_*^x)$ 
are finitely generated as modules over ${\bf R}$.}

{\bf   Theorem   4.} 
{\it If  
   ${\bf H}_0(x,\widehat{f}_x)    =
   {\bf R}[x]/(f(x))_x$ 
is a finitely generated module over ${\bf R}$, then\smallskip

\ \ \ $\exists  e(x_*,\widehat{f}_*^x) \in {\bf Z}(x_*,\widehat{f}_*^x): 
   \mathop{\top}\limits_{(y,\widehat{f}_y)}
   \det \left\| \begin{matrix}\ \nabla f(x,y) \ \\ 
   \ \widehat{f}_x-\widehat{f}_y \ \end{matrix} \right\|
   e(y_*,\widehat{f}_*^y)= 1+\partial \left[ t(x,\widehat{f}_x) \right] $
\smallskip \\
and 
%there 
it hold the statements of theorem 3.}\eject\noindent

{\footnotesize

\begin{enumerate}

\item {\it Seifullin, T. R.} 
Root functionals and root polynomials 
of a system of polynomials. (Russian)
Dopov. Nats. Akad. Nauk Ukra\"\i ni  -- 1995, -- no. 5, 5--8.
\item {\it Seifullin, T. R.} Root functionals and root relations 
of a system of polynomials. (Russian) 
Dopov. Nats. Akad. Nauk Ukra\"\i ni  -- 1995, -- no 6, 7--10.
\item {\it  Bourbaki, N.}  Alg\`ebre, Chapitre 10, Alg\`ebre homologique, Paris: Masson, 1980. 
\\
\end{enumerate}

\small{\noindent
{\it V. M. Glushkov Institute of Cybernetics of the NAS of Ukraine, Kiev
\hfill Received 29.10.96\medskip\\
E-mail: \ {\tt  timur\_sf@mail.ru}  
}

\end{document}